\documentclass[12pt]{amsart}
\usepackage{amscd,amssymb,times,fullpage}
\newtheorem{theorem}{Theorem}
\newtheorem{lemma}{Lemma}
\newtheorem{proposition}{Proposition}

\theoremstyle{definition}
\newtheorem{definition}{Definition}

\newtheorem{example}{Example}

\theoremstyle{remark}
\newtheorem{remark}{Remark}

\newcommand{\R}{\mathbb{ R}}

\newcommand{\GG}{\mathcal{ G}}

\newcommand{\FF}{\mathcal{ F}}
\newcommand{\RR}{\mathcal{ R}}

\newcommand\Symp{\operatorname{Symp}}

\newcommand\Ker{\operatorname{Ker}}

\newcommand\Diff{\operatorname{Diff}}

\catcode`\@=\active 

\begin{document}

\title{Stability theorems for symplectic and contact pairs}
\author{G.~Bande}
\address{Universit\`a degli Studi di Cagliari, Dip.~Mat., Via Ospedale 
72, 09129 Cagliari, Italy}
\email{gbande{\char'100}unica.it}
\author{P.~Ghiggini}
\address{Mathematisches Institut, Ludwig-Maximilians-Universit{\"a}t M{\"u}nchen,
Theresienstr.~39, 80333 M{\"u}nchen, Germany}
\email{ghiggini{\char'100}mathematik.uni-muenchen.de}
\author{D.~Kotschick}
\address{Mathematisches Institut, Ludwig-Maximilians-Universit{\"a}t M{\"u}nchen,
Theresienstr.~39, 80333 M{\"u}nchen, Germany}
\email{dieter{\char'100}member.ams.org}


\thanks{The authors are members of the {\sl European Differential
Geometry Endeavour} (EDGE), Research Training Network HPRN-CT-2000-00101,
supported by The European Human Potential Programme}
\date{\today; MSC 2000 classification: primary 53C15, 57R17, 57R30;
secondary 53C12, 53D35, 58A17}

\begin{abstract}
We prove Gray--Moser stability theorems for complementary pairs of 
forms of constant class defining symplectic pairs, contact-symplectic 
pairs and contact pairs. We also consider the case of 
contact-symplectic and contact-contact structures, in which the 
constant class condition on a one-form is replaced by the condition 
that its kernel hyperplane distribution have constant class in the 
sense of {\'E}.~Cartan.
\end{abstract}

\maketitle

\section{Introduction}

A symplectic pair on a smooth manifold $M$ is a pair of closed two-forms 
$\omega$ and $\eta$ of constant and complementary ranks, for which $\omega$ 
restricts as a symplectic form to the leaves of the kernel foliation of 
$\eta$, and vice versa. This definition from~\cite{KM}, see also~\cite{BK}, 
is analogous to that of contact pairs and of contact-symplectic pairs 
introduced in~\cite{Bande1,Bande2,BH}. In this paper we prove analogs of 
the stability theorems of Gray~\cite{Gray} and Moser~\cite{Moser} for 
all these structures.

Gray~\cite{Gray} proved that for a smoothly varying family $\xi_{t}$
of contact structures on a closed manifold there exists an isotopy
$\varphi_{t}$ with $\varphi_{t}^{*}\xi_{t}=\xi_{0}$. It is easy to see
that in general one cannot obtain an isotopy of contact forms, but
only of contact structures. Moser~\cite{Moser} proved a corresponding
theorem for families of volume forms or symplectic forms. In this
case one has to assume that the de Rham cohomology class of the forms
is constant, and one obtains an isotopy of forms.

In our situation of symplectic or contact pairs, we have to take into
account the diffeomorphism types of the foliations involved. It is
now well known that for geometric structures with underlying
foliations there can be no general stability theorems, as the
diffeomorphism type of the foliations may vary smoothly and
non-trivially, cf.~for example~\cite{Adachi,Golubev,MZ,Vogel}.
We shall give explicit examples in which such a variation happens for
symplectic pairs (with constant cohomology classes). The
constructions of~\cite{BK} then allow us to exhibit the same
phenomenon for contact-symplectic and for contact pairs.

To obtain stability theorems, we can only consider families in which the 
underlying foliations are constant. In this case we have to make assumptions 
not on the de Rham cohomology classes of forms, but on their cohomology classes 
in a refined or foliated cohomology. In a similar vein, Ghys~\cite{Ghys} has 
adapted Moser's argument for volume forms to leafwise volume forms for a
foliation $\FF$. To obtain the desired stability result he assumed
that the cohomology class of the volume forms is constant in the
leafwise cohomology of $\FF$. It is obvious that this adaptation also
works for leafwise symplectic forms whose cohomology classes are
constant in leafwise cohomology.

Our setup is slightly more complicated. We are not interested in isotopies 
of leafwise forms only, but in isotopies of forms on the whole manifold.
Assuming that the relevant cohomology classes are constant in de Rham
cohomology, it follows by restriction that they are constant in
leafwise cohomology. However, it will turn out that the desired
isotopy exists if and only if a stronger condition is satisfied: the
cohomology classes should be constant in the basic cohomology of the
complementary foliation. We shall give examples which show that this 
condition cannot be enforced by assuming only constancy in de Rham 
cohomology.

In Section~\ref{s:structures} we define contact-symplectic and 
contact-contact structures, which generalize contact-symplectic and 
contact pairs respectively by replacing the constant class conditions 
on one-forms by the constant class condition on their kernel hyperplane 
distributions. We give explicit examples of such structures which cannot 
be defined by contact-symplectic or contact pairs, and we prove the 
appropriate stability theorems. Note that each half of a contact-contact structure 
and the contact half of a contact-symplectic structure are hyperplane 
distributions of constant class in the sense of {\'E}.~Cartan. For such 
hyperplane distributions an analogue of Gray's theorem was proved by 
Montgomery and Zhitomirskii~\cite{MZ}, and our stability theorem for 
contact-contact structures reproves and generalizes their result.

\medskip
\noindent
{\bf Acknowledgement:} The third author is grateful to F.~Kamber for 
discussions in 1999 about adaptations of Moser's theorem for volume 
forms in the context of foliations. In particular, those discussions 
led to the conclusion that the basic cohomology of a foliation is the 
correct receptacle for the obstruction to isotopies of transverse volume 
forms, which foreshadowed the appearance of basic cohomology in this paper.

\section{Preliminaries on foliated cohomology}\label{s:prelim}

In this section we recall the notions from foliated cohomology that
will be needed for our arguments. We refer the reader to~\cite{Mo,S}
for further information.

Let $\FF$ be a smooth foliation on $M$. We denote by $I(\FF)$ the
(graded) ideal in the de Rham complex consisting of forms vanishing on $\FF$.
By the Frobenius theorem this ideal is closed with respect to the exterior
differential $d$. Thus we have a short exact sequence of differential complexes
\begin{equation}\label{eq:fol}
0\to I^*(\FF)\to \Omega^*(M)\to \Omega^*(M)/I^*(\FF)\to 0
\end{equation}
with an induced differential in the quotient complex. This quotient complex,
also denoted $\Omega^{*}(\FF)$, is the complex of leafwise forms for $\FF$.
Its cohomology, denoted $H^*(\FF)$,
is called the foliated, or leafwise, cohomology of $\FF$.

The exact sequence~\eqref{eq:fol} gives rise to a long exact sequence
in cohomology, whose connecting homomorphism is induced by the
exterior differential. In part, we have:
\begin{equation}\label{eq:les}
    \ldots\to H^{1}(M)\to H^{1}(\FF)\to
    H^{2}(I(\FF))\to H^{2}(M)\to
    H^{2}(\FF)\to\ldots
    \end{equation}
This exact sequence explains the relationship between different 
conditions one can impose on the cohomology classes of closed 
$2$-forms which are in the ideal of one foliation and are 
nondegenerate on another foliation, complementary to the first one.
For example, if $\FF$ is the kernel foliation of a closed two-form $\omega$ 
of constant rank, then $[\omega ]\in H^{2}(I(\FF))$ determines 
$[\omega ]\in H^{2}(M)$, but not the other way around.

Now suppose that $\FF$ has codimension $q$ and $\alpha$ is a global 
defining $q$-form. This is well-defined up to multiplication by a 
nowhere vanishing function. By the Frobenius theorem there exists a 
$1$-form $\beta$ such that $d\alpha= \alpha\land\beta$. Fixing $\alpha$, 
this $\beta$ is well-defined up to addition of an element $\tau\in I^1(\FF)$.
The following is well-known:
\begin{lemma}\label{Reebclass}
    The form $\beta$ is closed in $\Omega^*(M)/I^*(\FF)$ and its cohomology 
    class $[\beta]\in H^1(\FF)$, called the Reeb class of $\FF$, is 
    well-defined independent of choices.
\end{lemma}

The geometric meaning of the Reeb class is explained by the following 
result:
\begin{proposition} 
    The following conditions on a cooriented $p$-dimensional foliation $\FF$ are 
    equivalent:
\begin{enumerate}
\item There is a holonomy-invariant volume form on transversals.
\item The Reeb class $[\beta]\in H^1(\FF)$ vanishes.
\item $\FF$ can be defined by a closed form.
\end{enumerate}
\end{proposition}
Note that the foliations induced by a symplectic, contact or 
contact-symplectic pair are always defined by closed forms, and so 
have vanishing Reeb class.

We now introduce the complex of basic forms for a foliation $\FF$:
$$
\Omega^{*}_{b}(\FF)=\{\alpha\in\Omega^{*}(M) \ \vert \ i_{X}\alpha = 
0, \ L_{X}\alpha = 0 \ \forall X\in T\FF \} \ .
$$
This is a subcomplex of the de Rham complex, and its cohomology, 
denoted $H^{*}_{b}(\FF)$, is called the basic cohomology of $\FF$.
If $\FF$ is the kernel foliation of a closed $p$-form $\alpha$, then 
$\alpha$ is a basic form and defines a class $[\alpha ]\in H^{p}_{b}(\FF)$. 
In our stability theorems we shall use the second basic cohomology 
to formulate the cohomological condition for the existence of an 
isotopy of symplectic or contact-symplectic pairs with constant 
characteristic foliations. We shall use the following facts, see~\cite{S}:
\begin{example}\label{ex:vertical}
    If $\FF$ has as leaves the fibers of a fiber bundle $M\to 
    B$, then $H^{*}_{b}(\FF)$ is isomorphic to $H^{*}(B)$. In 
    particular it is finite-dimensional if the base $B$ is compact.
    \end{example}
\begin{example}\label{ex:horizontal}
    If $\FF$ is the horizontal foliation of a flat bundle with fiber 
    $F$ obtained by suspending $\rho\colon\pi_{1}(B)\to 
    \Diff(F)$, then $H^{*}_{b}(\FF)$ is isomorphic to the cohomology 
    of the subcomplex of $\rho$-invariant forms in $\Omega^{*}(F)$. In 
    particular, it is often infinite-dimensional.
    \end{example}

Note that the inclusion $i\colon\Omega^{*}_{b}(\FF)\to I(\FF)$ induces 
a map $i_{*}\colon H^{*}_{b}(\FF)\to H^{*}(I(\FF))$ in cohomology.
\begin{lemma}\label{l:incl}
    The map $i_{*}$ is an isomorphism in degree $1$, and is injective in degree $2$.
\end{lemma}
\begin{proof}
    A closed $1$-form in the ideal is basic, and so $i_{*}$ is 
    surjective in degree $1$.
    
    Suppose $\alpha$ is a basic $p$-form such that $i_{*}([\alpha])=0$. 
    Then there is a $(p-1)$-form in the ideal such that $d\beta =\alpha$. 
    For $p=1$ the form $\beta$ is a function and vanishes identically, so 
    $i_{*}$ is injective in degree $1$. For $p=2$ the form $\beta$ is a 
    one-form satisfying $\beta (X)=0$ for all $X\in T\FF$. Moreover, 
    $L_{X}\beta =i_{X}d\beta = i_{X}\alpha$, and this vanishes as $\alpha$ is 
    basic. Thus $\beta$ is basic and $i_{*}$ is injective in degree $2$ as well.
    \end{proof}
    
\section{Moser's theorem for symplectic pairs}

In this section we consider smoothly varying families of symplectic pairs
$(\omega_{t},\eta_{t})$ on a closed manifold $M$. This means that $\omega_{t}$
and $\eta_{t}$ are closed $2$-forms of constant and complementary ranks,
$\omega_{t}$ restricts as a symplectic form to the leaves of the kernel foliation
$\FF_{t}$ of $\eta_{t}$, and $\eta_{t}$ restricts as a symplectic form to the leaves
of the kernel foliation $\GG_{t}$ of $\omega_{t}$. In particular, $\FF_{t}$ and
$\GG_{t}$ are complementary smooth foliations.

We first show that there are such families with nondiffeomorphic
foliations.
\begin{example}\label{ex:fol}
    Let $(F,\omega_{F})$ be a closed symplectic manifold. Choose a nontrivial
    smoothly varying family $\varphi_{t}\in\Symp_{0}(F)$ of symplectomorphisms
    of $F$ isotopic to the identity. In fact, we may assume $\varphi_{0}=Id_{F}$.
    Let $M=F\times S^{1}$, endowed with the varying horizontal
    foliations given by considering $M$ as the mapping torus of
    $\varphi_{t}$. We can think of $N=M\times S^{1}$ as a foliated
    bundle over $T^{2}$.

    The form $\omega_{t}$ induced on $N$ from $\omega_{F}$ by suspending
    $\varphi_{t}$ and the pullback $\eta$ of a volume form on $T^{2}$ to $N$
    form a symplectic pair with constant $\eta$. The kernel foliation
    of $\omega_{t}$ however varies in a nontrivial way. For $t=0$ we
    assume that $\varphi_{0}$ is the identity, so all the leaves of
    the kernel foliation of $\omega_{t}$ are closed. For positive $t$
    we can introduce non-closed leaves by choosing the family
    $\varphi_{t}$ appropriately.
    \end{example}

    \begin{remark} \label{rem:BW}
    The above variation in the diffeomorphism type of the foliations can
    not be controlled by assumptions on the cohomology classes of the
    forms involved. Lemma~8 of~\cite{KM} shows that if the family
    $\varphi_{t}$ consists of Hamiltonian symplectomorphisms, then the
    cohomology class of $\omega_{t}$ is constant. Moreover, in this case it is
    integral if $\omega_{F}$ is integral on $F$.
    \end{remark}

From now on we only consider families of symplectic pairs for which
the foliations $\FF_{t}$ and $\GG_{t}$ are independent of $t$.
Thus $\omega_{t}$ and $\eta_{t}$ (or their maximal nonzero powers) are
defining forms for $\GG$ and $\FF$ respectively, for all $t$. To obtain 
an isotopy of symplectic pairs we have to assume at least that $[\omega_{t}]$ 
and $[\eta_{t}]$ are constant in the de Rham cohomology of $M$. The exact 
sequence~\eqref{eq:les} then shows that $[\omega_{t}]$ and $[\eta_{t}]$ are 
also constant in the leafwise cohomologies $H^{2}(\FF)$ and in $H^{2}(\GG)$ 
respectively. This is the analog of the assumption in Ghys's leafwise Moser 
argument~\cite{Ghys}. It turns out that we need more, and that the precise 
condition for the existence of an isotopy involves either the cohomology 
of the ideals of the foliations, or their basic cohomology. 

Note that $\eta_{t}$ is a basic form for $\FF$, and $\omega_{t}$ is 
basic for $\GG$. Lemma~\ref{l:incl} shows that assuming these forms 
to represent constant cohomology classes in the basic cohomology is in 
fact equivalent to the assumption that they are constant in the cohomology 
of the ideals. Thus, in the following stability theorem the assumption about 
the basic cohomology classes could be replaced by an equivalent assumption 
involving the cohomology of the ideals\footnote{{\it Mutatis 
mutandis} this remark applies to the stability theorems for 
contact-symplectic pairs and contact-symplectic structures as well.}.
\begin{theorem}\label{th:main}
    Let $(\omega_{t},\eta_{t})$ be a smooth family of symplectic pairs on 
    a closed smooth manifold $M$, such that the kernel foliations
    $\FF=\Ker(\eta_{t})$ and $\GG=\Ker(\omega_{t})$ are independent of
    $t\in [0,1]$. Then there exists an isotopy $\varphi_{t}$ with
    $\varphi_{t}^{*}\omega_{t}=\omega_{0}$ and $\varphi_{t}^{*}\eta_{t}=\eta_{0}$ 
    if and only if the basic cohomology classes $[\omega_{t}]\in H_{b}^{2}(\GG)$ 
    and $[\eta_{t}]\in H_{b}^{2}(\FF)$ are constant.
    \end{theorem}
The assumptions imply that the cohomology classes are constant in de
Rham cohomology. Clearly the isotopy preserves the foliations.
\begin{proof}
    Suppose that the isotopy $\varphi_{t}$ exists. Let $X_{t}$ be the
    time-dependent vector field generating it. Then we have
    $$
    0=\frac{d}{dt}(\varphi_{t}^{*}\omega_{t})=
    \varphi_{t}^{*}(\dot\omega_{t}+L_{X_{t}}\omega_{t})=
    \varphi_{t}^{*}(\dot\omega_{t}+di_{X_{t}}\omega_{t}) \ ,
    $$
    because $\omega_{t}$ is closed. It follows that
    $$
    \dot\omega_{t}=-di_{X_{t}}\omega_{t} \ .
    $$
    Let $Y\in T\GG$. Then 
    $i_{Y}(i_{X_{t}}\omega_{t})=-i_{X_{t}}(i_{Y}\omega_{t})$ vanishes 
    because $\omega_{t}$ is a basic form for $\GG$. Moreover, for 
    $Y\in T\GG$ we also have 
    $$
    L_{Y}(i_{X_{t}}\omega_{t})=i_{Y}di_{X_{t}}\omega_{t}=-i_{Y}\dot\omega_{t}=0
    $$
    because $\dot\omega_{t}$ is also basic for $\GG$. Thus we conclude that
    $[\dot\omega_{t}]=0\in H^{2}_{b}(\GG)$, and this condition is
    necessary for the existence of an isotopy. Similarly we derive the
    necessity of $[\dot\eta_{t}]=0\in H^{2}_{b}(\FF)$.

    Conversely, assume $[\omega_{t}]\in H^{2}_{b}(\GG)$ and $[\eta_{t}]\in
    H^{2}_{b}(\FF)$ are constant. Then there exist time-dependent one-forms 
    $\alpha_{t}\in\Omega^{1}_{b}(\GG)$ and $\beta_{t}\in\Omega^{1}_{b}(\FF)$,
    depending smoothly on $t$, such that $d\alpha_{t}=\dot\omega_{t}$ and 
    $d\beta_{t}=\dot\eta_{t}$. 
    Let the time-dependent vector field $X_{t}$ be defined by the equation
    \begin{equation}\label{eq:X}
    i_{X_{t}}(\omega_{t}+\eta_{t})=-\alpha_{t}-\beta_{t} \ .
    \end{equation}
    As $\omega_{t}+\eta_{t}$ is non-degenerate, $X_{t}$ exists and is
    uniquely determined. Rewriting~\eqref{eq:X} in the form
    \begin{equation}\label{eq:XX}
    i_{X_{t}}\omega_{t}+\alpha_{t} = -i_{X_{t}}\eta_{t}-\beta_{t} \ ,
    \end{equation}
    we have an identity of one-forms, where the left-hand side is in $I(\GG)$ and 
    the right-hand side is in $I(\FF)$. This means that both sides 
    vanish.

    Consider now the flow $\varphi_{t}$ of $X_{t}$. We have
    $$
    \frac{d}{dt}(\varphi_{t}^{*}\omega_{t})=
    \varphi_{t}^{*}(\dot\omega_{t}+L_{X_{t}}\omega_{t})=
    \varphi_{t}^{*}(\dot\omega_{t}+di_{X_{t}}\omega_{t}) =
    \varphi_{t}^{*}(\dot\omega_{t}-d\alpha_{t})=0
    $$
    by the vanishing of the left-hand side of~\eqref{eq:XX} and the 
    definition of $\alpha_{t}$. Thus 
    $\varphi_{t}^{*}\omega_{t}=\varphi_{0}$ for all $t\in [0,1]$.
    Similarly the vanishing of the right-hand side of~\eqref{eq:XX} 
    and the definition of $\beta_{t}$ show that 
    $\varphi_{t}^{*}\eta_{t}=\eta_{0}$, so that we have an
    isotopy between $(\omega_{1},\eta_{1})$ and $(\omega_{0},\eta_{0})$.
    \end{proof}

If we assume that $M$ is $4$-dimensional and consider symplectic pairs
$(\omega,\eta)$ with both two-forms of rank two, then these two-forms are
leafwise volume forms on $\FF$ and $\GG$ respectively. Thus, unlike for
general symplectic forms we can take convex combinations as follows. If
$(\omega_{1},\eta_{1})$ and $(\omega_{2},\eta_{2})$ are two symplectic
pairs with the same kernel foliations $\FF$ and $\GG$, and inducing the
same orientations on the leaves, then the combinations
$(t\omega_{1}+(1-t)\omega_{2},t\eta_{1}+(1-t)\eta_{2})$ also have these
properties for all $t\in [0,1]$. If we assume in addition that the two
pairs we start with represent the same cohomology class in
$H^{2}_{b}(\GG)\times H^{2}_{b}(\FF)$, then the argument above produces
an isotopy between them. Moreover, the argument works with arbitrary
parameters in a compact space, so that we obtain the following:
\begin{theorem}
    Let $(\omega,\eta)$ be a symplectic pair on a closed four-manifold.
    Denote by $\Diff_{+}(\FF,\GG,[\omega],[\eta])$ the group of
    diffeomorphisms preserving the oriented foliations $\FF$ and
    $\GG$ and the cohomology classes $[\omega]\in H^{2}_{b}(\GG)$
    and $[\eta]\in H^{2}_{b}(\FF)$. Denote by $\Symp(\omega,\eta)$ the
    group of diffeomorphisms preserving the forms $\omega$ and $\eta$.
    Then the inclusion
    $$
    \Symp(\omega,\eta)\hookrightarrow\Diff_{+}(\FF,\GG,[\omega],[\eta])
    $$
    is a weak homotopy equivalence for the $C^{\infty}$ topology on
    the two groups.
    \end{theorem}
Such a result was proved for volume forms by Moser~\cite{Moser}, and for
leafwise volume forms by Ghys~\cite{Ghys}.

Finally, we want to show that the assumptions about constancy of 
basic cohomology classes do not follow from the other assumptions in 
our setup. For this we use an example originally described by Ghys 
(unpublished) for different reasons.
\begin{example}
    Let $N$ be the $T^{2}$-bundle over $S^{1}$ with monodromy 
    $\phi =\begin{pmatrix} 1 & 1 \\ 0 & 1 \end{pmatrix}$. Then 
    $M=N\times S^{1}$ is a foliated $T^{2}$-bundle over $T^{2}$ with 
    area-preserving holonomy carrying an obvious symplectic pair 
    $(\omega,\eta)$, where $\omega$ is induced by a $\phi$-invariant area 
    form on the fiber, and $\eta$ is a pullback from the base. Let 
    $\GG=\Ker(\omega)$ be the horizontal foliation. By 
    Example~\ref{ex:horizontal} the basic cohomology $H^{2}_{b}(\GG)$ 
    is infinite-dimensional, and so we can certainly choose some 
    nontrivial $\dot\omega$ in the kernel of $H^{2}_{b}(\GG)\to 
    H^{2}(M)$. Setting $\omega_{t}=\omega +t\dot\omega$ for $\vert 
    t\vert$ small enough, we obtain a family of symplectic pairs 
    $(\omega_{t},\eta)$ which have fixed foliations and de Rham 
    cohomology classes, but which cannot be isotopic by 
    Theorem~\ref{th:main}. 
    \end{example}

\section{Gray's theorem for contact and contact-symplectic pairs}

Contact pairs and contact-symplectic pairs were introduced 
in~\cite{Bande1,Bande2,BH}, to which we refer for the basic properties.

We first recall the notion of class of a differential form.
\begin{definition}
The characteristic subspace of $\alpha$ at a point $p\in M$ is the subspace 
$C_{p}(\alpha)\subset T_{p}M$ given by the intersection of the 
kernels of $\alpha _p$ and of $d\alpha _p$. The class of $\alpha$ at $p$ is 
the codimension of $C_{p}(\alpha)$.
\end{definition}

A useful way for checking the class of a form of degree $1$ or $2$ is given 
by the following:
\begin{lemma}
The class of a closed $2$-form $\eta$ on $M$ is even, and equal to 
its rank. Thus the class of $\eta$ at $p$ is equal to $2k$ if and only if 
$\eta _{p}^{k}\neq 0$ and $\eta _{p}^{k+1}=0$.

A $1$-form $\alpha$ has class $2k+1$ at $p$ if and only if 
$\alpha _{p}\land (d\alpha_{p})^{k}\neq 0$ and $(d\alpha _{p})^{k+1}=0$. 
It has class $2k$ if and only if $(d\alpha _{p})^{k}\neq 0$ and 
$\alpha _{p}\land (d\alpha_{p})^{k}=0$.
\end{lemma}

It is easy to show that if $\omega$ has constant class on $M$ then its 
characteristic distribution is completely integrable. The resulting 
foliation is called the characteristic foliation of $\omega$. Its codimension is 
the class of $\omega$.

We consider manifolds equipped with a pair of forms of constant class 
satisfying certain additional conditions:  
\begin{definition}[\cite{Bande1,Bande2}]\label{d:cspair}
A pair $(\alpha, \eta )$, where $\alpha$ is a $1$-form and $\eta$ is a closed 
$2$-form, is called a contact-symplectic pair of type $(h,k)$ if the following 
conditions are satisfied: $\alpha\land (d\alpha)^{h}\land\eta ^{k}$ is a volume 
form, $(d\alpha)^{h+1}=0$ and $\eta ^{k+1}=0$.
\end{definition}
It follows from the definition that $\alpha$ has constant class $2h+1$ and $\eta$ 
has constant class $2k$, that their characteristic foliations are transverse, and 
that the leaves carry induced symplectic respectively contact structures. For these 
pairs there is a notion of Reeb vector fields, explained in the next proposition:
\begin{proposition}[\cite{Bande1,Bande2}]
For a contact-symplectic pair $(\alpha, \eta )$ there exists a unique vector field 
$R$, called the Reeb vector field of the pair, such that: $\alpha (R) 
=1$, $i_{R}d\alpha =0$ and $i_{R}\eta =0$.
\end{proposition}
Note that by definition $R$ is tangent to the kernel foliation of $\eta$.

We want to consider smoothly varying families of contact-symplectic
pairs $(\alpha_{t},\eta_{t})$. 
According to Definition~\ref{d:cspair} this means that $\alpha_{t}$ is a one-form 
of constant class $2k+1$, $\eta_{t}$ is a closed $2$-form of constant rank $2l$ 
such that the kernel foliations $\GG_{t}$ of $\alpha_{t}\land (d\alpha_{t})^{k}$ 
(which is equal to $\Ker(\alpha_t) \cap \Ker(d\alpha_t)$) and
$\FF_{t}$ of $\eta_{t}$ are complementary, $\alpha_{t}$ restricts as a
contact form to the leaves of $\FF_{t}$, and $\eta_{t}$ restricts as a symplectic
form to the leaves of $\GG_{t}$.

Performing a leaf-wise Boothby--Wang construction as in~\cite{BK} on the 
symplectic pairs in Example~\ref{ex:fol}, we obtain a smoothly varying family 
of contact-symplectic pairs with non-diffeomorphic foliations. Thus, to
obtain a stability result we need to assume again that the foliations $\FF_{t}$
and $\GG_{t}$ are independent of $t\in [0,1]$.

\begin{theorem}\label{th:main2}
    Let $(\alpha_{t},\eta_{t})$ be a family contact-symplectic pairs of type 
    $(h, k)$ on a closed smooth manifold $M$, such that the kernel foliations
    $\FF=\Ker(\eta_{t})$ and $\GG=\Ker(\alpha_{t}\land (d\alpha_{t})^{h})$ are 
    independent of $t\in [0,1]$. Then there exists an isotopy $\varphi_{t}$ with
    $\varphi_{t}^{*}\eta_{t}=\eta_{0}$ and $\varphi_{t}^{*}\alpha_{t}=f_{t}\alpha_{0}$ 
    for some smooth non-zero functions $f_{t}$ if and only if the basic cohomology
    class $[\eta_{t}]\in H^{2}_{b}(\FF)$ is constant.
    \end{theorem}
We do not give the proof here, because Theorem~\ref{th:main2} is a special case 
of Theorem~\ref{th:main4} proved below in the more general context of 
contact-symplectic structures.

Next we recall the definition of contact pairs:
\begin{definition}[\cite{Bande1,BH}]\label{d:cpair}
A pair $(\alpha, \beta)$ of $1$-forms on a manifold is said to be a contact pair 
of type $(h,k)$ if the following conditions are satisfied: $\alpha\land 
(d\alpha)^{h}\land\beta\land (d\beta)^{k}$ is a volume form, $(d\alpha)^{h+1}=0$
and $(d\beta)^{k+1}=0$.
\end{definition}
The forms $\alpha$ and $\beta$ have constant class $2h+1$ and $2k+1$ respectively, 
and the leaves of their characteristic foliations have induced contact structures.

Manifolds with contact pairs carry a pair of Reeb vector fields:
\begin{proposition}[\cite{Bande1,BH}]\label{p:Reeb}
For a contact pair $(\alpha,\beta)$ there exist two commuting vector fields $A$, 
$B$ uniquely determined by the following  conditions: $\alpha(A)=\beta(B)=1$,
$\alpha(B)=\beta(A)=0$ and $i_{A}d\alpha=i_{A}d\beta=i_{B}d\alpha=i_{B}d\beta=0$.
\end{proposition}

Again performing leaf-wise Boothby--Wang constructions as in~\cite{BK} on the 
symplectic pairs in Example~\ref{ex:fol}, we obtain smoothly varying families 
of contact pairs with non-diffeomorphic foliations. Thus, to obtain a stability 
result we need to assume that the foliations $\FF$ and $\GG$ are independent of 
$t\in [0,1]$.
\begin{theorem}\label{th:main3}
    Let $(\alpha_{t},\beta_{t})$ be a family contact pairs of type $(h, k)$ on a 
    closed smooth manifold $M$, such that the kernel foliations $\FF=
    \Ker(\alpha_{t}\land (d\alpha_{t})^{h})$ and $\GG=\Ker(\beta_{t}\land 
    (d\beta_{t})^{k})$ are independent of $t\in [0,1]$. Then there exists an 
    isotopy $\varphi_{t}$ with $\varphi_{t}^{*}\alpha_{t}=f_t \alpha_{0}$ and
    $\varphi_{t}^{*}\beta_{t}=g_t \beta_{0}$.
    \end{theorem}
This is a special case of Theorem~\ref{th:main5}, which we prove 
below in the more general context of contact-contact structures.

\section{Contact-symplectic and contact-contact structures}\label{s:structures}

In this section we introduce contact-symplectic and contact-contact structures, 
which are generalizations of contact-symplectic and of contact pairs 
respectively, and we prove stability theorems for these new structures.
\begin{definition}
A pair $(\xi, \eta )$ of a hyperplane field $\xi$ and a closed 2-form $\eta$ on 
a $(2n+1)$-dimensional manifold $M$ is called a contact-symplectic structure of 
type $(h, n-h)$ if $\xi$ can be defined by a $1$-form $\alpha$ such that the 
following conditions are satisfied: $\alpha\land (d\alpha)^{h}\land\eta^{n-h}
\neq 0$, $\alpha\land (d\alpha)^{h+1}=0$ and $\eta^{n-h+1}=0$.
\end{definition}
Clearly these conditions do not depend on the choice of defining form $\alpha$.
Since $d \alpha + \eta$ has constant rank $2n$ and its kernel is not 
contained in the kernel of $\alpha$, the following definition makes 
sense:
\begin{definition}
The Reeb vector field $R$ of $(\alpha, \eta )$ is the unique vector field such 
that $\alpha (R)=1$ and $i_R (d\alpha+ \eta)=0$.
\end{definition}
The Reeb vector field is not an invariant of the contact-symplectic structure,
but depends on the choice of form $\alpha$.

It is well-known that given a hyperplane field $\xi=\Ker \alpha$ of 
constant class $h$, i.~e.~with $\alpha\land (d\alpha)^{h+1}=0$ and 
$\alpha\land (d\alpha)^{h} \neq 0$ for one and therefore all defining 
forms $\alpha$, the distribution $\Ker (\alpha\land (d\alpha)^{h})$ is 
completely integrable, and locally there exists a form $\beta$ of constant class 
(as a form) such that $\Ker (\alpha\land (d\alpha)^{h})=\Ker (\beta)$.
Thus a contact-symplectic structure, like a contact-symplectic pair, determines 
two complementary foliations which have symplectic and contact leaves respectively.
It is also clear that a contact-symplectic pair determines a contact-symplectic 
structure, but the converse is not true, as shown by the following example:
\begin{example}\label{excsstruct}
Let $M$ be any closed oriented three-manifold which does not fiber over the circle, 
endowed with a taut foliation $\FF$ given by the kernel of a one-form $\alpha$. 
There exists a one-form $\beta$ such that $d\alpha=\alpha \land \beta$, 
and $\beta$ represents the Reeb class of $\FF$. Since we have chosen $M$ not to 
fiber over the circle, $\alpha$ can not be closed, and the Reeb class 
$[\beta] \in H^1(\FF)$ is non-zero.

As we have chosen $\FF$ to be taut, there exists a closed two-form $\eta$ such 
that $\alpha\land\eta >0$. It follows that $(\Ker(\alpha),\eta)$ is 
a contact-symplectic structure of type $(0,1)$ which cannot be defined by a 
contact-symplectic pair of forms of constant class (as forms).
\end{example}

\begin{example}
    For a slightly more complicated example take $M$ as above and 
consider the manifold $N=M \times T^{2}$ endowed with the pair 
$(\Ker(\gamma), \eta)$, where $\gamma=\cos x \alpha + \sin x dy$ and $x$ and $y$ are 
standard coordinates on the torus. An easy calculation shows that 
this is a contact-symplectic structure of type 
$(1,1)$, and that $d(\gamma\land d\gamma)=\gamma\land d\gamma 
\land (2\cos^2 x \cdot \beta)$, where $\beta$ represents the Reeb class of 
$\FF$. This means that the foliation defined by $\gamma \land d\gamma$ has 
non-zero Reeb class $2 \cos^2 x \cdot \beta$ and therefore cannot be defined 
by a closed form. This implies that the contact-symplectic structure 
does not come from an underlying contact-symplectic pair.
\end{example}

Next we generalize the notion of contact pairs:
\begin{definition}
A pair $(\xi, \sigma)$ of hyperplane fields on a $(2n+2)$-dimensional manifold $M$ 
is called a contact-contact structure of type $(h, n-h)$ if there are 
defining one-forms $\alpha$ and $\beta$ for which the following conditions 
are satisfied: $\alpha\land (d\alpha)^{h}\land\beta\land 
(d\beta)^{n-h}\neq 0$, $\alpha\land (d\alpha)^{h+1}=0$ and 
$\beta\land (d\beta)^{n-h+1}=0$.
\end{definition}
The conditions in the definition are independent of the choice of 
defining forms $\alpha$ and $\beta$.
We define a Reeb distribution as follows:
\begin{definition} 
    The Reeb distribution $\RR$ consists of the tangent vectors $Y$ 
    satisfying the equation $(i_Y(d\alpha +d\beta))\vert_{\xi\cap\sigma}=0$.
\end{definition}
It is easy to see that this is a smooth distribution of rank two, 
which depends on the choice of $\alpha$ and $\beta$. We can unravel 
the definition as follows. The $2$-form $d\alpha+d\beta$ has rank at 
least $2n$, because it is non-degenerate on the transverse 
intersection of the hyperplanes $\xi$ and $\sigma$. If its rank at a 
point is exactly $2n$, then at that point the fiber of the Reeb 
distribution $\RR$ is the kernel of $d\alpha+d\beta$. If the rank of 
$d\alpha+d\beta$ at a point is not $2n$, then the form is symplectic 
in an open neighbourhood of that point in $M$, and on that 
neighbourhood $\RR$ is the symplectic orthogonal of $\xi\cap\sigma$.  
\begin{definition}\label{reebdist}
The Reeb vector fields $A$, $B$ of $(\alpha, \beta )$ are the unique 
vector fields tangent to the Reeb distribution $\RR$ such that 
$\alpha (A)=\beta (B)=1$, $\alpha (B)=\beta (A)=0$.
\end{definition}

The contact-contact structures with commuting Reeb vector fields are 
special:
\begin{proposition}
    The following conditions on defining forms $\alpha$ and $\beta$ 
    for a contact-contact structure $(\xi,\sigma)$ are equivalent:
    \begin{itemize}
        \item[(i)] The Reeb vector fields $A$ and $B$ commute.
        \item[(ii)] The form $d\alpha+d\beta$ has constant rank $2n$.
        \end{itemize}
These conditions imply
    \begin{itemize}
        \item[(iii)] The Reeb distribution $\RR$ is integrable.
        \end{itemize}
    \end{proposition}
\begin{proof}
    It is clear that (i) implies (iii), so we only have to prove the 
    equivalence of (i) and (ii).
    
    By definition, we have $d\alpha (A,B)=-\alpha ([A,B])$ and 
    $d\beta (A,B)=-\beta ([A,B])$. 
    
    If $A$ and $B$ commute, then 
    $d\alpha+d\beta$ vanishes on $\RR$, the span of $A$ and $B$. If 
    $d\alpha+d\beta$ had maximal rank $2n+2$, then it would have to be 
    nondegenerate on $\RR$, the symplectic orthogonal of $\xi\cap\sigma$. 
    
    Conversely, if $d\alpha+d\beta$ has rank $2n$ everywhere then it 
    defines a $2$-dimensional foliation, whose leaves are tangent to 
    $\RR$. Thus $i_{A}(d\alpha+d\beta)$ and $i_{B}(d\alpha+d\beta)$ 
    vanish, so that by the formulae above $[A,B]=0$.
    \end{proof}
The Reeb distribution and the Reeb vector fields are not invariants of the
contact-contact structure, but they depend on the forms $(\alpha, \beta)$.
If the contact-contact structure can be defined using a contact pair 
$(\alpha, \beta)$ of constant class in the sense of Definition~\ref{d:cpair}, 
then the Reeb vector fields defined above coincide with those in 
Proposition~\ref{p:Reeb}, they commute, and the Reeb distribution $\RR$ is integrable.

\begin{example}
    Consider again a three-manifold $M$ with a taut foliation $\FF$ that cannot be 
    defined by a closed form, as in Example~\ref{excsstruct}. The two-form $\eta$ 
    which is positive on the leaves of the foliation can be taken to represent an 
    integral cohomology class. Let $\pi\colon N\to M$ be the circle bundle 
    with Euler class $[\eta]$, endowed with a connection one-form $\gamma$ with 
    $d\gamma=\pi^{*}\eta$. If $\alpha$ is any defining form for $\FF$, then 
    $\pi^{*}\alpha$ and $\gamma$ are defining forms for a contact-contact structure 
    on $N$ which cannot be defined by a contact pair. Nevertheless, the Reeb vector 
    fields for $\pi^{*}\alpha$ and $\gamma$ commute.
    \end{example}

For both contact-symplectic and contact-contact structures the Reeb classes 
in the sense of Lemma~\ref{Reebclass} of the characteristic foliations are 
the obstructions to finding a contact-symplectic or contact pair which defines 
the structure.
    
To prove the stability theorems, we need the following:
\begin{lemma}\label{lemma1}
Let $\alpha$ be a one-form such that $\alpha\land (d\alpha)^{h} \neq 0$ and 
$\alpha\land (d\alpha)^{h+1} = 0$ and let $X$ be a vector field tangent to 
the kernel of $\alpha$. Let $\FF$ be the kernel foliation of 
$\alpha\land (d\alpha)^{h}$. Then $i_ X d\alpha$ is in the ideal $I(\FF)$.
\end{lemma}
\begin{proof}
Contracting $\alpha\land d\alpha^{h+1} = 0$ with $X$ we obtain 
$(i_X d\alpha)\land\alpha\land d\alpha^{h}=0$, which simply means that 
$i_ X d\alpha$ is in the ideal $I(\FF)$.
\end{proof}

Given a contact-contact  structure as above, we can also define leafwise Reeb
vector fields as the vector fields tangent to the characteristic 
foliations, which on each leaf are the Reeb vector field of the 
induced contact structure:
\begin{definition}
The leafwise Reeb vector fields $R_{\alpha}$ and $R_{\beta}$ are the unique vector 
fields tangent to $\Ker(\beta\wedge (d\beta)^{n-h}$ and $\Ker(\alpha\wedge (d\alpha)^{h}$ 
respectively satisfying the conditions:
$$
\alpha (R_{\alpha})=\beta (R_{\beta})=1 \ ,
$$
$$
(i_{R_{\alpha}}d\alpha) \wedge \beta \wedge (d\beta)^{n-h} =0 \ ,
$$
$$
(i_{R_{\beta}}d\beta) \wedge \alpha \wedge (d\alpha)^{h} =0 \ .
$$
\end{definition}
The following proposition gives a relation between the Reeb vector 
fields in the sense of Definition~\ref{reebdist} and the leafwise Reeb vector fields:
\begin{proposition}
The leafwise Reeb vector field $R_{\alpha}$ is the projection of $A$ 
to $\Ker(\beta\wedge (d\beta)^{n-h})$ along $\Ker(\alpha\wedge (d\alpha)^{h})$, 
and similarly for $R_{\beta}$.
\end{proposition}
\begin{proof}
We show that $A$ projects to $R_{\alpha}$, the other case being completely symmetric.
The definition of $A$ implies that there exists a function $f\colon M \to \R$ such 
that $i_A (d \alpha + d \beta)= f \beta$. Wedging both sides of this equality 
with $\beta\wedge (d \beta)^{n-k}$ we obtain $(i_A d \alpha 
)\wedge\beta\wedge (d \beta)^{n-k} = 0$ because $(i_A d\beta 
)\wedge\beta\wedge (d \beta)^{n-k} = 0$ by Lemma~\ref{lemma1}. Thus the vector field 
$X=A-R_{\alpha}$ satisfies $(i_X d \alpha )\wedge\beta\wedge (d 
\beta)^{n-h}=0$. It also satisfies $(i_X d\alpha)\wedge\alpha\wedge (d\alpha)^h$ by 
Lemma~\ref{lemma1} because $i_X \alpha =0$. Thus we conclude $i_X d \alpha =0$. The 
vector field $X$ is therefore tangent to $\Ker(\alpha\wedge (d\alpha)^{h})$ because 
$(i_X \alpha)\wedge\alpha\wedge (d \alpha)^h=0$. 
\end{proof}
In general the Reeb vector fields $A$ and $B$ and the leafwise Reeb 
vector fields $R_{\alpha}$ and $R_{\beta}$ are different. We 
introduced $A$ and $B$ because they will be useful in the proofs of 
the stability theorems.

Finally, here is the stability theorem for contact-symplectic structures:
\begin{theorem}\label{th:main4}
    Let $(\Ker(\alpha_{t}),\eta_{t})$ be a smooth family of contact-symplectic 
    structures of type $(h, k)$ on a closed smooth manifold $M$, such that the 
    kernel foliations $\FF=\Ker(\eta_{t})$ and 
    $\GG=\Ker(\alpha_{t}\land (d\alpha_{t})^{h})$ are independent of $t\in [0,1]$. 
    Then there exists an isotopy $\varphi_{t}$ such that 
    $\varphi_{t}^{*}\eta_{t}=\eta_{0}$ and $\varphi_{t}^{*}\alpha_{t}=f_{t}\alpha_{0}$ 
    for some smooth non-zero functions $f_{t}$, if and only if the basic cohomology 
    class $[\eta_{t}]\in H^{2}_{b}(\FF)$ is constant.
    \end{theorem}
\begin{proof}
    The proof that constancy of the cohomology class is necessary is the same as in 
    the proof of Theorem~\ref{th:main}.
    
    Conversely, assume $[\eta_{t}]\in H^{2}_{b}(\FF)$ is constant. This implies that
    there is a time-dependent one-form $\beta_{t}$ which is basic 
    for $\FF$ and such that $d\beta_{t}=\dot\eta_{t}$. Define a time-dependent vector
    field $X_{t}$ by the equations
    \begin{equation}\label{eq:X4}
    i_{X_{t}}(d\alpha_{t}+\eta_{t})=-\beta_{t}+ \mu _{t} \cdot \alpha _{t} - \dot\alpha _{t}  \ ,
    \end{equation}
    \begin{equation}\label{ex:XX4}
    i_{X_{t}}(\alpha_{t})=0 \ ,
    \end{equation}
where $\mu _{t}=\dot\alpha _{t} (R_t) + \beta_{t} (R_t)$, for $R_{t}$ 
the Reeb vector field.
To show that this definition makes sense, consider for 
   every $t$ and every point $p\in M$ the following linear map given by 
    $d\alpha_{t}+\eta_{t}$:
\begin{alignat*}{1}
    \Ker(\alpha_{t})_{p} &\longrightarrow (R_{t})_{p}^{\perp}\\
v &\longmapsto i_{v}(d\alpha_{t}+\eta_{t})_p \ ,
\end{alignat*}
where $(R_{t})_{p}^{\perp}\subset T^{*}_{p}M$ is the annihilator of the value 
of the Reeb vector field $(R_{t})_{p}$. Because the dimensions of these vector 
spaces are the same, and $d\alpha_{t}+\eta_{t}$ is nondegenerate on 
$\Ker(\alpha_{t})_p$, this map is an isomorphism. Thus an $X_{t}$ as above exists and 
is uniquely determined because the right-hand side of~\eqref{eq:X4} 
vanishes on the Reeb vector field.
    
From~\eqref{eq:X4} we obtain:
    $$
    i_{X_{t}}d\alpha_{t}-\mu _{t} \cdot \alpha _{t} + \dot\alpha _{t}=
    -i_{X_{t}}\eta_{t}-\beta_{t} \ .
    $$
    As $\alpha_{t}$ and $\dot\alpha_{t}$ are in $I(\GG)$, $i_{X_{t}}d\alpha_{t}$ is 
    in $I(\GG)$ by Lemma~\ref{lemma1}, and $i_{X_{t}}\eta_{t}$ and $\beta_{t}$ are 
    in $I(\FF)$, we conclude that each side of the previous equation vanishes.
 Consider now the flow $\varphi_{t}$ of $X_{t}$. We have
    $$
    \frac{d}{dt}(\varphi_{t}^{*}\eta_{t})=
    \varphi_{t}^{*}(\dot\eta_{t}+L_{X_{t}}\eta_{t})=
    \varphi_{t}^{*}(\dot\eta_{t}+di_{X_{t}}\eta_{t}) =
    \varphi_{t}^{*}(\dot\eta_{t}-d\beta_{t})=0
    $$
    by the definition of $\beta_{t}$, and
$$
\frac{d}{dt}(\varphi_{t}^{*}\alpha_{t})=\varphi_{t}^{*}(\dot\alpha_{t}+L_{X_{t}}\alpha_{t})=
\varphi_{t}^{*}(\dot\alpha_{t}+i_{X_{t}}d\alpha_{t})=\varphi_{t}^{*}(\mu 
_{t}\cdot\alpha _{t}) \ .
$$
This implies 
$$
\varphi_{t}^{*}\alpha_{t}=f_t \alpha_0
$$
with $f_t=\exp (\int_{0} ^{t}(\varphi_{s} ^{*}\mu _{s})ds)$ and
$$
\varphi_{t}^{*}\eta_{t}=\eta_0 .
$$
This completes the proof.
\end{proof}

Now we consider the case of contact-contact structures:
\begin{theorem}\label{th:main5}
    Let $(\Ker(\alpha_{t}),\Ker(\beta_{t}))$ be a family contact-contact structures 
    of type $(h,n-h)$ on a closed smooth manifold $M$, such that the kernel foliations 
    $\FF=\Ker(\alpha_{t}\land (d\alpha_{t})^{h})$ and $\GG=\Ker(\beta_{t}\land 
    (d\beta_{t})^{k})$ are independent of $t\in [0,1]$. Then there exists an 
    isotopy $\varphi_{t}$ so that $\varphi_{t}^{*}\alpha_{t}=f_t \alpha_{0}$ and
    $\varphi_{t}^{*}\beta_{t}=g_t \beta_{0}$.
    \end{theorem}
\begin{proof}
Let $A_t$, $B_t$ the Reeb vector fields of the pair, as in Definition~\ref{reebdist}. 
Define a vector field $W_t$ by the following equations:
$$
    i_{W_{t}}(d\alpha_{t}+d\beta_{t})=\mu _{t} \cdot \alpha _{t} - 
    \dot\alpha _{t} + \nu _{t} \cdot \beta _{t} - 
    \dot\beta _{t}\ ,
$$
$$
    \alpha_{t}(W_t)=0 \ ,
$$
$$
    \beta_{t}(W_t)=0 \ ,
$$
where $\mu _{t}=\dot\alpha_{t} (A_t) + \dot\beta_{t} (A_t)$ and 
$\nu _{t}=\dot\alpha_{t} (B_t) + \dot\beta_{t} (B_t)$.

To see that $W_t$ is well-defined consider the map
\begin{alignat*}{1}
    \Ker(\alpha_{t})_{p}\cap \Ker(\beta_{t})_{p} &\longrightarrow 
    (A_{t})_{p}^{\perp}\cap (B_{t})_{p}^{\perp}\\
v &\longmapsto i_{v}(d\alpha_{t}+d\beta_{t})_p \ .
\end{alignat*}
Now $d\alpha_{t}+d\beta_{t}$ is nondegenerate on $\Ker(\alpha_{t})\cap\Ker(\beta_{t})$, 
and so this map is an isomorphism. The choice of $\mu_{t}$ and $\nu _{t}$ implies that 
$\mu_{t}\cdot\alpha_{t} -\dot\alpha_{t}+\nu _{t}\cdot\beta_{t}-\dot\beta_{t}$ vanishes 
on $A_{t}$ and on $B_{t}$. Thus $W_t$ exists and is unique.

From the first equation we obtain:
$$
    i_{W_{t}} (d\alpha_{t})- \mu _{t} \cdot \alpha _{t} + 
    \dot\alpha _{t} = - i_{W_{t}} d\beta_{t} + \nu _{t} \cdot \beta _{t} - 
    \dot\beta _{t}\ .
$$
Here each side vanishes because on the left-hand side all forms are in $I(\FF)$ 
and on the right-hand side all are in $I(\GG)$.
Let $\varphi_{t}$ be the flow generated by $W_{t}$. 
Then after differentiation by $t$ we obtain:
$$
\varphi_{t}^{*}(\alpha _{t})=f_t\alpha_{0} \ ,
$$
$$
\varphi_{t}^{*}(\beta  _{t})=g_t\beta_{0}  \ ,
$$
where $f_t=\exp (\int_{0} ^{t}(\varphi_{s}^{*}\mu _{s})ds)$ and 
$g_t=\exp (\int_{0} ^{t}(\varphi_{s}^{*}\nu _{s})ds)$.
\end{proof}

\bigskip

\bibliographystyle{amsplain}


\end{document}